\documentclass[12pt]{article}
\usepackage{amssymb}
\usepackage{graphicx}
\usepackage{amsfonts}
\usepackage{latexsym}
\usepackage{amsthm}
\usepackage{amsmath}

\usepackage{amssymb, amscd}

\usepackage{url}

\newtheorem{problem}{Problem}
\newtheorem{theo}[problem]{Theorem}
\newtheorem*{rem}{Remark}

\newtheorem{defin}[problem]{Definition}
\newtheorem{prop}[problem]{Proposition}
\newtheorem{cor}[problem]{Corollary}
\newtheorem{lema}[problem]{Lemma}
\newtheorem{exam}[problem]{Example}

\begin{document}
\date{March 20, 2013}
 \title{{Rotation number of a unimodular cycle:\\
an elementary approach}}

\author{{Rade  T.\ \v Zivaljevi\' c}\\ {\small Mathematical Institute}\\[-2mm] {\small SASA, Belgrade}
\\[-2mm]{\small rade$@$mi.sanu.ac.rs}}

\maketitle
\begin{abstract}
We give an elementary proof of a formula expressing the rotation
number of a cyclic unimodular sequence $L=u_1u_2\ldots u_d$ of
lattice vectors $u_i\in \mathbb{Z}^2$ in terms of arithmetically
defined local quantities. The formula has been originally derived
by A.~Higashitani and M.~Masuda (arXiv:1204.0088v2 [math.CO]) with
the aid of the Riemann-Roch formula applied in the context of
toric topology. These authors also demonstrated that a generalized
version of the `Twelve-point theorem' and a generalized Pick's
formula are among the consequences or relatives of their result.
Our approach emphasizes the role of `discrete curvature
invariants' $\mu(a,b,c)$, where $\{a,b\}$ and $\{b,c\}$ are bases
of $\mathbb{Z}^2$,  as fundamental discrete invariants of {\em
modular lattice geometry}.
\end{abstract}

\section{Introduction}\label{sec:introduction}

The following theorem of A.~Higashitani and M.~Masuda, proved in
\cite{H-M}, is a close relative of the remarkable `Twelve-point
theorem' \cite{P-R}, \cite{crs}, \cite{Fulton}, \cite{H-S},
\cite{K-N}. Like its predecessors, the `Twelve-point theorem' and
Pick's formula, it is an intriguing and easily formulated
statement about a sequence of lattice vectors, their rotation
(winding) number and the associated, arithmetically defined local
quantities.

\begin{theo}\label{thm:glavna_teorema}
The rotation number $Rot(L)$ of a cyclic unimodular sequence $L =
u_1u_2\ldots u_d$ can be calculated as the weighted sum
\begin{equation}\label{eqn:glavna-formula}
Rot(L) = \frac{1}{12}\mu(L) + \frac{1}{4}\nu(L) =
\frac{1}{12}\sum_{i=1}^d\, \mu(u_{i-1},u_i,u_{i+1}) +
\frac{1}{4}\sum_{i=1}^d\, \nu(u_i,u_{i+1})
\end{equation}
of locally defined quantities $\mu(L)$ and $\nu(L)$ where
$\nu(u_i,u_{i+1}):={\rm det}(u_i,u_{i+1})\in\{-1,+1\}$ and
$\mu(u_{i-1},u_i,u_{i+1})=a_i\in \mathbb{Z}$ is the integer
determined by the equation
\[
{\rm det}(u_{i-1},u_{i})u_{i-1} + {\rm det}(u_{i},u_{i+1})u_{i+1}
+ a_iu_i = 0.
\]
\end{theo}

\medskip
Theorem~\ref{thm:glavna_teorema} may appear at first sight as
quite elementary and not difficult to comprehend, however it has a
deeper meaning and significance. Like its relative (and a
consequence) the `Twelve-point theorem', it is situated at the
crossroads of several mathematical areas, illuminating and
offering new perspectives on `well understood' mathematical
concepts.

The first proof \cite{H-M} (see also \cite[Section 5]{Ma}) of the
formula (\ref{eqn:glavna-formula}) was based on a Riemann-Roch
type theorem (Noether formula) where the integers $a_i$ appeared
as the self-intersection numbers of the corresponding homology
classes of the associated `topological toric manifold'.

\medskip
Our first objective in this paper is to give a conceptual and
elementary proof of Theorem~\ref{thm:glavna_teorema} which is
based on a systematic analysis of the invariant $\mu(a,b,c)$. The
second objective is to pave the way for the hypothetical higher
dimensional analogues of Theorem~\ref{thm:glavna_teorema}. For
this reason the exposition emphasizes the study of invariants
$\mu(a,b,c)$ and their higher dimensional versions $\mu_j(a,
\mathfrak{b}, a')$ (Section~\ref{sec:appendix}) as `discrete
curvature invariants' situated within {\em unimodular lattice
geometry}. This point of view is similar to the approach of
O.~Karpenkov  to  `lattice trigonometry', as part of his
investigation of `lattice geometry invariants' \cite{K-1, K-2}.

\medskip
A different, short and elementary computational proof of
Theorem~\ref{thm:glavna_teorema} was subsequently included in the
new version of the paper \cite{H-M}. Their Lemma~1.3.\  fits in
nicely in our approach so we took the opportunity to shorten the
original proof  (\cite{Z}, arXiv:1209.4981v1[math.CO]), retaining
its transparency and conceptuality.

We observe in passing (Section~\ref{sec:other-proofs}) that the
proof of Poonen and Rodrigu\-ez-Villegas of the `Twelve-point
theorem' \cite{P-R}, based on the properties of the holomorphic
function (modular form) $\Delta(z)$, can also be used as the basis
of a proof of Theorem~\ref{thm:glavna_teorema}. The variety of
proofs and methods applied are perhaps an indication that this
result deserves a further exploration so the paper ends with some
open questions and conjectures.

\section{Introductory definitions and remarks}
\label{sec:remarks-definitions}

\subsection{Unimodular sequences}
\label{sec:definicije}

A sequence $L=u_1u_2\ldots u_d$ of lattice vectors is called {\em
unimodular} if $\{u_i,u_{i+1}\}$ is a basis of the lattice
$\mathbb{Z}^2$ for each $i$ or equivalently if ${\rm
det}(u_i,u_{i+1})\in \{-1,+1\}$ for each $i=1,\ldots,d-1$.
Geometrically this condition means that for each $i$ the triangle
$Ou_iu_{i+1}$ does not contain lattice points aside from the
vertices.

\medskip
A unimodular sequence $L=u_1u_2\ldots u_d$ is called {\em cyclic}
if ${\rm det}(u_d,u_1)\in\{-1,+1\}$. A cyclic unimodular sequence
$L$ (of length $d$) naturally defines a $d$-periodic unimodular
sequence $W = \ldots LLL\ldots = \ldots u_{-1}u_0u_1u_2\ldots
u_du_{d+1}\ldots$ where $u_i = u_j$ if $i \equiv j\, \, \mbox{\rm
(mod\, $d$)}$.

\subsection{Local and global $\mu$-invariants}
\label{sec:mu-and-nu}

The invariants $\mu(L)$ and $\nu(L)$ of a cyclic unimodular
sequence $L = u_1\ldots u_d$ are already introduced in
Theorem~\ref{thm:glavna_teorema} as the sums
\begin{equation}\label{eqn:mu-nu-definicija}
\mu(L) = \sum_{i=1}^d\, \mu(u_{i-1},u_i,u_{i+1}) \qquad \nu(L) =
\sum_{i=1}^d\, \nu(u_i,u_{i+1}).
\end{equation}
The $\mu$-invariant of a unimodular sequence $(u,v,w)$ is
described as the unique integer $a=\mu(u,v,w)$ determined by the
equation
\begin{equation}\label{eqn:mu-invarijanta-definicija}
{\rm det}(u,v)u + {\rm det}(v,w)w + av = 0.
\end{equation}
Together with the associated $\nu$-invariant $\nu(u,v):={\rm
det}(u,v)$ this is a basic discrete angle invariant of (planar)
{\em modular lattice geometry}. Higher dimensional analogues of
these invariants are introduced and discussed in
Section~\ref{sec:appendix} (see Definition~\ref{def:gen-mu}).

A possible ambiguity arises if $L=u_1u_2u_3$ is a cyclic
unimodular sequence. For this reason the term `$\mu$-invariant' is
reserved for the number $\mu(u_1,u_2,u_3)$ (local $\mu$-invariant)
while $\mu(L)=\mu(u_1,u_2,u_3) + \mu(u_2,u_3,u_1) +
\mu(u_3,u_1,u_2)$ is the corresponding global $\mu$-invariant.

\subsection{Rotation number}
\label{sec:rotation-number}

Let $P = P(a_1,\ldots, a_d)$ be a closed, oriented, polygonal
curve in the plane with points $a_i$ as vertices and
$\overline{a_i,a_{i+1}}=[a_i, a_{i+1}]$ as oriented edges
$(a_{d+1}:= a_1)$. If the origin $O$ is not on $P$ it has a {\em
rotation number}\ (or winding number) defined by,
\begin{equation}\label{eqn:rotation-number-discrete}
{ Rot}(P) = \frac{1}{2\pi}\sum_{i=1}^d
\nu(a_i,a_{i+1})\angle(a_iOa_{i+1})
\end{equation}
where $\nu(a_i,a_{i+1})$ is the sign of the determinant ${\rm
det}(a_i,a_{i+1})$ and $\angle(a_iOa_{i+1})$ is the measure of the
angle $a_iOa_{i+1}$.

\medskip
Given a cyclic, unimodular sequence $L=u_1u_2\ldots u_d$ the
associated {\em closed unimodular polygon}\ is
$P_L=P(u_1,u_2,\ldots, u_d)$.  The {\em rotation number}\, ${
Rot}(L)$ of $L$ is by definition the rotation number of the
polygonal curve $P_L$.

\medskip
It is well known that ${ Rot}(P)$ can be defined with the aid of
elementary homology theory. We do not need this definition here
but we shall occasionally use the term {\em unimodular cycle}
$[L]$ to describe the collection
$\{\overline{u_iu_{i+1}}\}_{i=1}^d$ of oriented edges of $P_L$
which may be written also as a formal sum,
\[ [L] :=
\overrightarrow{u_1u_2} + \overrightarrow{u_2u_3} + \ldots +
\overrightarrow{u_{d-1}u_d} + \overrightarrow{u_du_1}.
\]
In this context the decomposition $[L] = [L_1] + [L_2] $, that
appears in Section~\ref{sec:additivity}, simply indicates that
$[L] = [L_1]\cup [L_2]$ is a union of sets with signed elements
(multisets) where the elements with differen sign, i.e.\ the edges
with different orientation, are supposed to cancel out. The
following proposition is, in light of the equation
(\ref{eqn:rotation-number-discrete}), an immediate consequence.

\begin{prop}\label{prop:rot-additivity} If $L, L_1$ and $L_2$ are cyclic,
unimodular sequences such that $[L] = [L_1] + [L_2]$ then ${
Rot}(L)={ Rot}(L_1) + { Rot}(L_2)$.
\end{prop}

\subsection{Modular vs. distance geometry}
\label{sec:mod-vs-distance}

The formula (\ref{eqn:glavna-formula}) exhibits interesting
similarities and differences when compared with more conventional
formulas for the rotation number of a planar curve $\Gamma$, say
the formula (\ref{eqn:rotation-number-discrete}) or its smooth
analogue,
\begin{equation}\label{eqn:arktan}
Rot(\Gamma) = \frac{1}{2\pi}\int_{\Gamma}\frac{xdy -
ydx}{x^2+y^2}.
\end{equation}
All these formulas are local, expressing the rotation number as a
sum (integral) of quantities defined locally on the curve. All of
them apply to closed curves, considering that each cyclic
unimodular sequence $L$ is associated a closed unimodular
polygonal line $P_L$. On closer inspection they even display an
interesting analogy between the set $S^1$ of unit vectors in
$\mathbb{R}^2$ and the set $Prim(\mathbb{Z}^2)$ of primitive
lattice vectors. Indeed, the local quantity integrated in
(\ref{eqn:arktan}) is the (infinitesimal) angle $d\theta =
d\arctan(y/x)$ where $\theta :
\mathbb{R}^2\setminus\{0\}\rightarrow S^1$ is the function
$\theta(v)= \|v\|^{-1}v$, whereas a similar role of ``discrete
angle functions'' is in formula (\ref{eqn:glavna-formula}) played
by the functions $\mu$ and $\nu$.

It is fairly certain that formulas (\ref{eqn:glavna-formula}) and
(\ref{eqn:arktan}) (or (\ref{eqn:glavna-formula}) and
(\ref{eqn:rotation-number-discrete})) cannot be directly related,
however it is likely that both can be derived from more general,
unifying principles. Discovering these principles appears to be an
interesting research problem.

\section{Properties of the invariant $\mu$}
\label{sec:mu}

\begin{prop}\label{prop:mu-invariant-formula}
Suppose that $(u,v,w)$ is an ordered triple of lattice vectors
such that $(u,v)$ and $(v,w)$ are ordered bases of $\mathbb{Z}^2$.
Then,
\begin{equation}\label{eqn:mu-formula}
\mu(u,v,w) = {\rm det}(u,v){\rm det}(v,w){\rm det}(w,u).
\end{equation}
\end{prop}

\medskip\noindent
{\bf Proof:} By definition $\mu(u,v,w)=a$ is the integer
determined by the equation
\[
{\rm det}(u,v)u + {\rm det}(v,w)w + a v = 0.
\]
A multiplication of both sides of this equality by ${\rm
det}(u,v){\rm det}(v,w)$ and a comparison with the equation
(\ref{eqn:3-vektora}) in the following lemma
(Lemma~\ref{lema:tri-vektora}) yields the desired formula. \hfill
$\square$

\begin{lema}\label{lema:tri-vektora} Suppose that $u,v,w\in \mathbb{R}^2$. Then,
\begin{equation}\label{eqn:3-vektora}
{\rm det}(u,v)w + {\rm det}(v,w)u + {\rm det}(w,u)v = 0
\end{equation}
Moreover, if ${\rm Span}(u,v,w)= \mathbb{R}^2$ then {\rm
(\ref{eqn:3-vektora})} is up to a scalar factor the only linear
dependence between vectors $u, v, w$.
\end{lema}

\begin{figure}[hbt]
\centering
\includegraphics[scale=0.40]{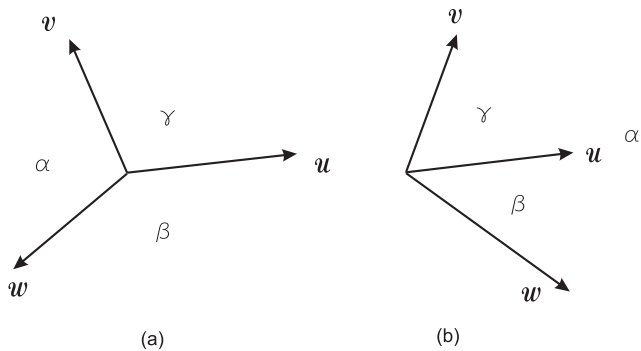}
\caption{Cyclic unimodular sequences of length $3$.}
\label{fig:3-vektora}
\end{figure}

\begin{cor}\label{cor:isometry}
If $u,v,w$ is a unimodular sequence then $\mu(w,v,u)=-\mu(u,v,w)$
and $\mu(-u,-v,-w) = \mu(u,v,w)$. More generally $\mu(O(u), O(v),
O(w))={\rm det}(O)\mu(u,v,w)$ for each $O\in GL(\mathbb{Z},2)$.
\end{cor}

As an illustration let us check Theorem~\ref{thm:glavna_teorema}
for the case of cyclic unimodular sequences of length $3$ and
special sequences of length $4$.

\begin{exam}\label{exam:3-vektora}
{\rm Assume that $L=(u,v,w)$ is a cyclic unimodular sequence. Let
$\alpha := {\rm det}(v,w), \beta := {\rm det}(w,u)$ and $\gamma :=
{\rm det}(u,v)$ (Figure~\ref{fig:3-vektora}). Then $\mu(L) =
3\alpha\beta\gamma$ (by
Proposition~\ref{prop:mu-invariant-formula}) and $\nu(L) = \alpha
+ \beta + \gamma$. Assume initially that $0\in {\rm
conv}\{u,v,w\}$ and that the sequence is positively oriented
(Figure~\ref{fig:3-vektora}~(a)) which means that
$\alpha=\beta=\gamma=1$. Then ${Rot}(u,v,w)=1$ while
$\mu(L)=\nu(L)=3$. Similarly if $-\alpha=\beta=\gamma=1$
(Figure~\ref{fig:3-vektora}~(b)) then $Rot(L)=0, \mu(L)=-3$ and
$\nu(L)=1$. This is sufficient to establish the formula,
$$Rot(L) =
({1}/{12})\mu(L) + ({1}/{4})\nu(L) = (1/4)(\alpha\beta\gamma +
\alpha + \beta + \gamma)$$ since both sides of the second equality
are symmetric with respect to $\alpha, \beta, \gamma$ and both
change the sign if the orientation of $L$ is reversed.  }
\end{exam}

\begin{figure}[hbt]
\centering
\includegraphics[scale=0.40]{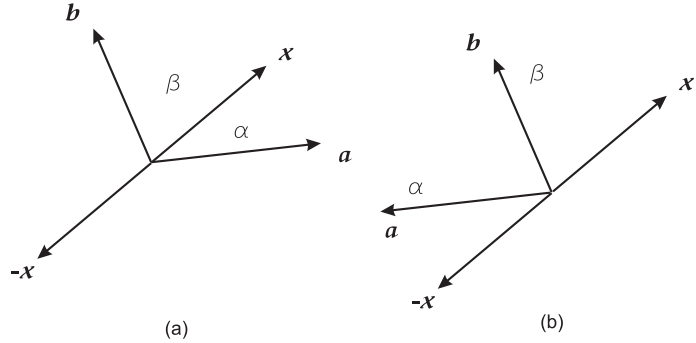}
\caption{Special unimodular sequences of length $4$.}
\label{fig:4-vektora}
\end{figure}

\begin{exam}\label{exam:4-vektora}
{\rm Let $L = \{a,x,b,-x\}$ be a special cyclic, unimodular
sequence of length $4$ (Figure~\ref{fig:4-vektora}). Observe that
$\mu(L)=0$ since for example $\mu(b,-x,a)=-\mu(a,x,b)$. If $\alpha
= {\rm det}[a,x]$ and $\beta = {\rm det}[x,b]$ then
$\nu(L)=2(\alpha + \beta)$ and the equality ${
Rot}(L)=(1/2)(\alpha + \beta) = (1/4)\nu(L)$ follows by inspection
of Figure~\ref{fig:4-vektora}. }
\end{exam}

\begin{rem}
Suppose $a,x,b$ is a unimodular sequence. Then  $a\pm x,x,b,
a,x,b\pm x$ and $a,-x,b$ are also unimodular and,
\[
\begin{array}{ccc}
\mu(a\pm x,x,b) & = & \mu(a,x,b) \mp {\rm det}(x,b) \\
\mu(a,x,b\pm x) & = & \mu(a,x,b) \mp {\rm det}(a,x)
\end{array}
\]
\end{rem}

\subsection{Exchange lemma and the Jacobi identity}
\label{sec:exchange}

\begin{prop}{\rm (Exchange Lemma)}\label{prop:first-exchange}
Suppose that $a,x,b$ and $u,x,v$ are (not necessarily cyclic)
unimodular sequences. Then the sequences $a,x,v$ and $u,x,b$ are
also unimodular and
\begin{equation}\label{eqn:first-exchange}
\mu(a,x,b)+\mu(u,x,v) = \mu(a,x,v) + \mu(u,x,b).
\end{equation}
\end{prop}

\medskip\noindent
{\bf Proof:} Let us write the defining equations of  the
corresponding $\mu$-invariants.
\[
\begin{array}{lcc}
{\rm det}(a,x)a + {\rm det}(x,b)b + \mu(a,x,b)x & = & 0 \\
{\rm det}(u,x)u + {\rm det}(x,v)v + \mu(u,x,v)x & = & 0 \\
{\rm det}(a,x)a + {\rm det}(x,v)v + \mu(a,x,v)x & = & 0 \\
{\rm det}(u,x)u + {\rm det}(x,b)b + \mu(u,x,b)x & = & 0
\end{array}
\]
The equality (\ref{eqn:first-exchange}) is obtained by subtracting
the third and fourth equation from the sum of the first two
equations. \hfill $\square$

\begin{prop}{\rm (Jacobi identity)}\label{prop:Jacobi}
Suppose that $(a,x), (b,x), (c,x)$ are bases of the lattice
$\mathbb{Z}^2$ which implies the unimodularity of $C = (a,x,b), A
= (b,x,c)$ and $B = (c,x,a)$. Then,
\begin{equation}\label{eqn:mu-3}
\mu(a,x,b) + \mu(b,x,c) + \mu(c,x,a) = 0
\end{equation}
\begin{equation}\label{eqn:mu-3-bis}
\mu(a,x,b) = \mu(a,x,c) + \mu(c,x,b)
\end{equation}
\end{prop}

\medskip\noindent
{\bf Proof:} Since $\mu(z,y,x)= -\mu(x,y,z)$
(Corollary~\ref{cor:isometry}) the equalities (\ref{eqn:mu-3}) and
(\ref{eqn:mu-3-bis}) are equivalent so it is sufficient to prove
only one of them.

The Pl\" ucker relation associated to the $2\times 4$ matrix  $M =
[a\, b\, c\, x]$ is
\[
{\rm det}(a,b){\rm det}(c,x) - {\rm det}(a,c){\rm det}(b,x) + {\rm
det}(a,x){\rm det}(b,c) = 0.
\]
On multiplying both sides of this equality by ${\rm det}(a,x){\rm
det}(b,x){\rm det}(c,x)$ we obtain on the left-hand side the
expression
\[
{\rm det}(a,b){\rm det}(a,x){\rm det}(b,x) - {\rm det}(a,c){\rm
det}(a,x){\rm det}(c,x) + {\rm det}(b,c){\rm det}(b,x){\rm
det}(c,x)
\]
which by Proposition~\ref{prop:mu-invariant-formula} immediately
leads to the equation (\ref{eqn:mu-3}). \hfill $\square$

\section{The additivity of $\mu$ and $\nu$  }
\label{sec:additivity}

In this section we describe a reduction procedure for cyclic
unimodular sequences which serves as a basis for an inductive
proof of the formula (\ref{eqn:glavna-formula}). The general idea
is to express a given sequence $L$, or rather the associated cycle
$[L]$ (Section~\ref{sec:rotation-number}), as a sum $L = L_1 +
L_2$ of simpler sequences such that $\mu(L) = \mu(L_1) + \mu(L_2)$
and $\nu(L_1)+\nu(L_2)$. Since the rotation number has a similar
behavior $Rot(L)= Rot(L_1)+Rot(L_2)$
(Proposition~\ref{prop:rot-additivity}), this eventually reduces
(\ref{eqn:glavna-formula}) to the case of cyclic, unimodular
sequences of length $2$ and $3$ where the formula is easily
verified (Example~\ref{exam:3-vektora}).

\subsection{Removing self-intersections}\label{sec:selfintersections}

\begin{figure}[hbt]
\centering
\includegraphics[scale=0.65]{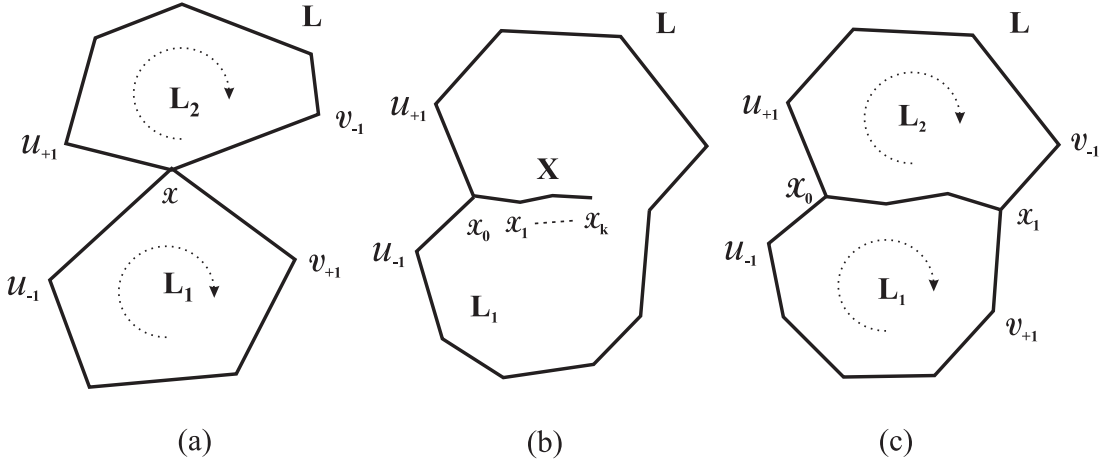}
\caption{Splitting and pruning of unimodular cycles.}
\label{fig:samopresek-1}
\end{figure}

Suppose that $L = \,.\,.\, u_{-1}xu_{+1}\,.\,.\,
v_{-1}xv_{+1}\,.\,.\,$ is a cyclic unimodular sequence which has a
{\em self-intersection} in the sense that a primitive vector $x$
appears twice in $L$, Figure~\ref{fig:samopresek-1}~(a). Then $L$
allows a ``shortcut decomposition'' $L = L_1 + L_2$ (more
precisely $[L] = [L_1] + [L_2]$) where $L_1 = \,.\,.\,
u_{-1}xv_{+1}\,.\,.\,$ and $L_2 = \,.\,.\, v_{-1}xu_{+1}\,.\,.\, $
are subsequences of $L$, Figure~\ref{fig:samopresek-1}~(a).

\begin{prop}\label{prop:selfintersection}
Suppose that a cyclic unimodular sequence $L$ has a
self-intersection (Figure~\ref{fig:samopresek-1}~(a)). If $L =
L_1+L_2$ is the associated shortcut decomposition then $\mu(L) =
\mu(L_1)+\mu(L_2)$ and $\nu(L)=\nu(L_1)+\nu(L_2)$.
\end{prop}

\medskip\noindent
{\bf Proof:} The relation for $\nu$ is obvious while the
additivity of $\mu$ follows from the Exchange Lemma
(Proposition~\ref{prop:first-exchange}). Indeed, the difference
$\mu(L) - \mu(L_1)-\mu(L_2)$ reduces to,
\[
\mu(u_{-1},x,u_{+1}) + \mu(v_{-1},x,v_{+1}) - \mu(u_{-1},x,v_{+1})
- \mu(v_{-1},x,u_{+1}) = 0.
\]

\begin{prop}\label{prop:back-and-forth}
Suppose that $L = \,.\,.\, u_{-1}x_0x_1\,.\,.\, x_k\,.\,.\,
x_1x_0u_{+1}\,.\,.\,$ is a cyclic unimodular sequence which has a
fragment $X$ where it goes back and forth
(Figure~\ref{fig:samopresek-1}~(b)). Let $L_1=\,.\,.\,
u_{-1}x_0u_{+1}\,.\,.\,$ be the cyclic unimodular sequence
obtained from $L$ by removing the segment $X = x_0x_1\,.\,.\,
x_k\,.\,.\, x_1x_0$. Then,
\[
\mu(L) = \mu(L_1)  \qquad \mbox{\rm\em and} \qquad  \nu(L) =
\nu(L_1).
\]
\end{prop}

\medskip\noindent
{\bf Proof:} The second equality is obvious while the first
follows from the antisymmetry of the $\mu$-invariant
(Corollary~\ref{cor:isometry}) and the Jacobi identity
(Proposition~\ref{prop:Jacobi}). Indeed,
\[
\mu(L) - \mu(L_1) = \mu(u_{-1},x_0,x_{1}) + \mu(x_1,x_0,u_{+1}) -
\mu(u_{-1},x_0,u_{+1}) = 0
\]
by an application of identity (\ref{eqn:mu-3-bis}). \hfill
$\square$

\subsection{Removing triangles}\label{sec:removing-triangles}

Suppose that a cyclic unimodular sequence $L = \,.\,.\,
v_{-1}v_0v_{+1}\,.\,.\,$ has a built-in unimodular triangle $L_1 =
v_{-1}v_0v_{+1}$; this situation arises if ${\rm
det}(v_{-1},v_{+1})\in \{-1,+1\}$. Then the sequence $L_0$
obtained from $L$ by deleting the vector $v_0$ is also a cyclic
unimodular sequence and there is a homological decomposition
$L=L_0+L_1$ (meaning that $[L]=[L_0]+[L_1]$). The following
proposition shows that the invariants $\mu$ and $\nu$ behave in
the expected way.

\begin{figure}[hbt]
\centering
\includegraphics[scale=0.45]{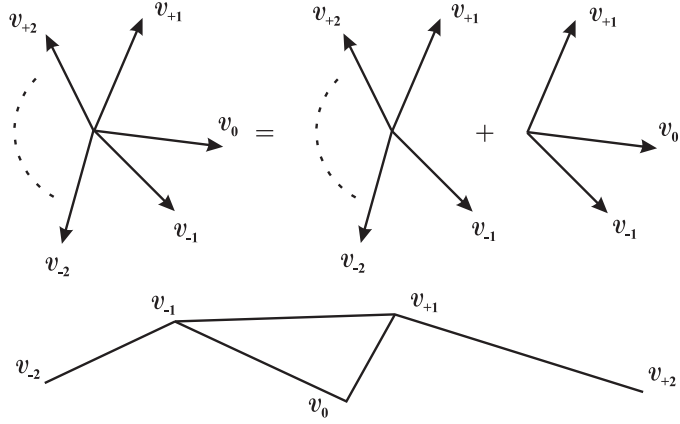}
\caption{Removing a triangle.} \label{fig:5-vektora}
\end{figure}

\begin{prop}\label{prop:redukcija-triangle}
Let us assume that a cyclic unimodular sequence $L$ has three
consecutive vectors $v_{-1}, v_0, v_{+1}$ such that $v_{-1},
v_{+1}$ is a basis of the lattice $\mathbb{Z}^2$. Then $L_0$
obtained from $L$ by deleting the vector $v_0$ and the three-term
sequence $L_1 = (v_{-1}, v_0, v_{+1})$ are both cyclic unimodular
sequences and,
\begin{equation}
\mu(L) = \mu(L_0) + \mu(L_1) \qquad \nu(L) = \nu(L_0) + \nu(L_1)
\end{equation}
\end{prop}

\medskip\noindent
{\bf Proof:} As before the relation $\nu(L)=\nu(L_0)+\nu(L_1)$ is
straightforward. The difference $\mu(L)-\mu(L_0)-\mu(L_1)$ is
equal to $A-B-C$ where
\[
A = \mu(v_{-2},v_{-1},v_{0}) + \mu(v_{-1},v_0,v_{+1}) +
\mu(v_{0},v_{+1},v_{+2})
\]
\[
B =  \mu(v_{-2},v_{-1},v_{+1}) + \mu(v_{-1},v_{+1},v_{+2})
\]
\[
C = \mu(v_{-1},v_{0},v_{+1}) + \mu(v_{0},v_{+1},v_{-1}) +
\mu(v_{+1},v_{-1},v_{0}).
\]
All terms in $A-B-C$ cancel out as a consequence of the Jacobi
identity (Proposition~\ref{prop:Jacobi}, equation
(\ref{eqn:mu-3-bis})) applied on vertices $v_{-1}$ and $v_{+1}$ of
the triangle $L_1 = v_{-1}v_0v_{+1}$ (Figure~\ref{fig:5-vektora}).
\hfill $\square$

\subsection{Splitting cycles}\label{sec:splitting-cycles}

Proposition~\ref{prop:redukcija-triangle} is a special case of a
general ``unimodular cycle splitting principle'' which is also a
consequence of the Jacobi identity and the Exchange Lemma.

\begin{prop}\label{prop:splitting-cycles}
Suppose that a cyclic unimodular sequence $$L = \,.\,.\,
u_{-1}x_0u_{+1}\,.\,.\, v_{-1}x_1v_{+1}\,.\,.\,$$ (as depicted in
Figure~\ref{fig:samopresek-1}~(c)) admits a unimodular shortcut $X
= x_0\,.\,.\, x_1$ giving rise to two new cyclic unimodular
sequences $L_1 = \,.\,.\, u_{-1}Xv_{+1}\,.\,.\,$ and
$L_2=\,.\,.\,v_{-1}X'u_{+1}\,.\,.\,$ (where $X':= x_1\,.\,.\,
x_0$). Then,
$$
\mu(L) = \mu(L_1) + \mu(L_2) \qquad {\mbox{\rm\em and}} \qquad
\nu(L) = \nu(L_1) + \nu(L_2).
$$
\end{prop}

\medskip\noindent
{\bf Proof:} The second part of the proposition (referring to the
invariant $\nu$) is again obvious. The first equality can be
proved directly, along the lines of the proof of
Proposition~\ref{prop:redukcija-triangle}, or deduced as a
consequence of Propositions~\ref{prop:selfintersection} and
\ref{prop:back-and-forth}.

Indeed, if $X = x_0Yx_1$ (and $X' = x_1Y'x_0$) then there exist
decompositions $L \cong Ax_0Bx_1, L_1 \cong AX, L_2 \cong BX'$,
where $\cong$ expresses equality of cyclic unimodular sequences
`up to a cyclic permutation'. Then by
Proposition~\ref{prop:back-and-forth} $\mu(L) =
\mu(Ax_0Yx_1Y'x_0Bx_1)$ and by
Proposition~\ref{prop:selfintersection} $\mu(Ax_0Yx_1Y'x_0Bx_1) =
\mu(L_1) + \mu(L_2)$. \hfill$\square$

\section{The proof of Theorem~\ref{thm:glavna_teorema}}
\label{sec:proofs}

Before commencing the proof of Theorem~\ref{thm:glavna_teorema} we
formulate a lemma which was used in \cite{H-M} for a similar
purpose, see also \cite[Section~2.5]{Fulton} (Exercises) where a
lemma of this kind was used as a combinatorial basis for
classification of smooth toric surfaces.

\begin{lema}{\em (\cite[Lemma~1.3]{H-M})}\label{lema:za-dokaz}
If $L = u_1u_2\ldots u_d$ is a cyclic, unimodular sequence of
length $d\geq 3$ then there exist three consecutive vectors
$u_{j-1}, u_{j}, u_{j+1}$ such that $$\mu(u_{j-1}, u_{j},
u_{j+1})\in\{0,+1, -1\}.$$
\end{lema}

\medskip\noindent
{\bf Proof of Theorem~\ref{thm:glavna_teorema}}: The proof is by
induction on the length $d$ of the cyclic, unimodular sequence $L
= u_1u_2\ldots u_d$. In the rather trivial case $d=2$,
$\mu(L)=\nu(L)={ Rot}(L)=0$. The case $d=3$ is established by
Example~\ref{exam:3-vektora} so let us suppose that $d\geq 4$. By
Lemma~\ref{lema:za-dokaz} we know that there exists an index $j$
such that $\mu(u_{j-1}, u_j, u_{j+1}) = a_j\in \{-1,0,+1\}$.

Assume that $a_j\in \{-1,+1\}$. Then ${\rm
det}(u_{j-1},u_{j+1})\in\{-1,+1\}$ as well (see
(\ref{eqn:mu-formula}) in
Proposition~\ref{prop:mu-invariant-formula}) and we are allowed to
use Proposition~\ref{prop:redukcija-triangle}. More explicitly,
$L$ has a decomposition $L = L_1 + L_2$ (in the sense of
Section~\ref{sec:additivity}) where $L_1 = u_1\ldots
u_{j-1}u_{j+1}\ldots u_d$ and $L_2$ is a triangular, cyclic,
unimodular sequence $L_2 = u_{j-1}u_ju_{j+1}$. By the induction
hypothesis the formula (\ref{eqn:glavna-formula}) holds for both
$L_1$ and $L_2$ hence, in light of the additivity of the
invariants $\mu, \nu$ and ${ Rot}$
(Propositions~\ref{prop:rot-additivity} and
\ref{prop:redukcija-triangle}), it holds also for $L$.

\begin{figure}[hbt]
\centering
\includegraphics[scale=0.40]{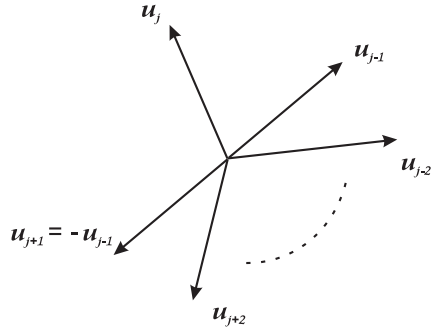}
\caption{The case $a_j=0$.} \label{fig:zero-case}
\end{figure}
In the case $a_j=0$ we again have a (one-step) unimodular
shortcut, this time from $u_{j-1}$ to $u_{j+2}$
(Figure~\ref{fig:zero-case}). We again have a decomposition $L =
L_1 + L_2$ where $L_1 = u_1\ldots u_{j-1}u_{j+2}\ldots u_d$ and
$L_2$ is a special cyclic, unimodular sequence $L_2 =
u_{j-1}u_ju_{j+1}u_{j+2}$ of length $4$. We know by
Example~\ref{exam:4-vektora} that formula
(\ref{eqn:glavna-formula}) holds for $L_2$ hence, again by
additivity (Propositions~\ref{prop:rot-additivity} and
\ref{prop:splitting-cycles}) and the induction hypothesis, it
holds also for $L$. \hfill $\square$

\subsection{Other proofs of Theorem~\ref{thm:glavna_teorema}}
\label{sec:other-proofs}

Here we briefly outline  the history of
Theorem~\ref{thm:glavna_teorema} and some of its immediate
predecessors. We record some of the main ideas used for its proof,
illustrating the diversity of methods and the relevance of this
result for different mathematical disciplines.

\medskip
The original proof of A.~Higashitani and M.~Masuda \cite[v2]{H-M},
see also \cite[Section~5]{Ma}, relied on {\em toric topology}
(Noether formula). This elegant proof nicely illustrates the
ongoing project of modifying or replacing some of the methods of
algebraic geometry (toric varieties) by the ideas of equivariant
algebraic topology (topological toric manifolds).

Probably the earliest appearance of a result that can be directly
linked to Theorem~\ref{thm:glavna_teorema} is
\cite[Section~2.5]{Fulton}. Fulton does not formulate it as a
separate statement about winding numbers, however in one of the
exercises he states that the equality $a_1+\ldots + a_d = 3d-12$
follows form the assumption that  $L = u_1u_2\ldots u_d$ is (in
our notation) a cyclic unimodular sequence such that $\nu(u_i,
u_{i+1})=+1$ and $\mu(u_{i-1}, u_i, u_{i+1})=a_i$ for each $i$.
His proof is inductive and uses a `topological constraint' on $L$
similar to Lemma~\ref{lema:za-dokaz}.

B.~Poonen and F.~Rodriguez-Villegas in their very interesting
paper \cite{P-R} use the properties of the modular form
$\Delta(z)$ to construct a homomorphism $\Phi :
\widetilde{SL}_2(\mathbb{Z})\rightarrow \mathbb{Z}$, where
$\widetilde{SL}_2(\mathbb{Z})$ is an extension of the group
${SL}_2(\mathbb{Z})$. They formulate a generalization of the
`Twelve-point theorem' to so called `legal loops' (cyclic
unimodular sequences) but they stop short of formulating
Theorem~\ref{thm:glavna_teorema}. However, in the course of the
proof they demonstrate how $\Phi$ can be used to evaluate the
associated rotation number. They also observe (Section~10) that
the $\mu$-invariant can be interpreted as a `combinatorial
analogue of an exterior angle' and discuss a possible connection
with the Gauss-Bonnet theorem.

The second proof of A.~Higashitani and M.~Masuda \cite[v3]{H-M} is
closer in spirit to our proof in Section~\ref{sec:proofs}. This
proof is also inductive, short and elementary, and may be a useful
alternative for those interested in a direct, computational proof
of Theorem~\ref{thm:glavna_teorema}.

\section{Generalizations and questions}\label{sec:appendix}

It is quite natural  to ask for higher dimensional analogues of
the formula (\ref{eqn:glavna-formula}) and other results from
previous sections. As a step in this direction we introduce the
concept of a {\em unimodular map} and discuss higher dimensional
generalizations of invariants  $\mu$ and $\nu$ from
Sections~\ref{sec:remarks-definitions} and \ref{sec:mu}.

\begin{defin}\label{def:unimodular-map}
Suppose that $M^d$ is an oriented, triangulated, closed manifold.
A map $\phi : M^d \rightarrow \mathbb{R}^{d+1}$ is called {\em
unimodular} if for each $d$-simplex $(a_0,a_1,\ldots, a_d)$ in
$M^d$ the set $\{\phi(a_0), \phi(a_1),\ldots, \phi(a_d)\}$ is a
basis of the lattice $\mathbb{Z}^{d+1}\subset \mathbb{R}^{d+1}$.
\end{defin}

\medskip\noindent
{\bf Problem~1:} Suppose that $M^d$ is an oriented, triangulated,
closed manifold. Let $\phi : M^d \rightarrow \mathbb{R}^{d+1}$ be
a {\em unimodular} map. If $[M]\in H_d(M^d; \mathbb{Z})$ is the
fundamental class of $M$ then
$$\phi_\ast([M])\in H_d(\mathbb{R}^{d+1}\setminus\{0\};
\mathbb{Z})\cong \mathbb{Z}$$ and $\phi_\ast([M])= k \iota$ for
some integer $k$ where $\iota$ is a generator of
$H_d(\mathbb{R}^{d+1}\setminus\{0\}; \mathbb{Z})$). The problem is
to find a formula expressing the integer $k$ in terms of locally
defined quantities, analogous to (\ref{eqn:glavna-formula}).

\medskip
There are natural candidates for the `locally defined quantities'
which are relatives and higher dimensional analogues of the
invariant $\mu$ introduced in Section~\ref{sec:mu}. Suppose that
$\mathcal{A} =\{a_1, a_2, b_1,\ldots, b_n\}$ is a collection of
integer vectors such that both $\mathcal{A}\setminus\{a_1\}$ and
$\mathcal{A}\setminus\{a_2\}$ are bases of $\mathbb{Z}^n$. Let
\begin{equation}\label{eqn:bases}
\nu_1 = [a_1; \mathfrak{b}] := {\rm det}(a_1,b_1,\ldots, b_n),\,
\nu_2 = [a_2; \mathfrak{b}] := {\rm det}(a_2,b_1,\ldots, b_n).
\end{equation}
It is not difficult to show (generalizing the formula
(\ref{eqn:mu-invarijanta-definicija}) from
Section~\ref{sec:mu-and-nu})  that there exists a unique sequence
of integers $\mu_1, \ldots, \mu_n$ such that,
\begin{equation}\label{eqn:Oda}
[a_1; \mathfrak{b}]a_1 - [a_2; \mathfrak{b}]a_2 +\mu_1b_1 + \ldots
+ \mu_nb_n =0.
\end{equation}
\begin{figure}[hbt]
\centering
\includegraphics[scale=0.45]{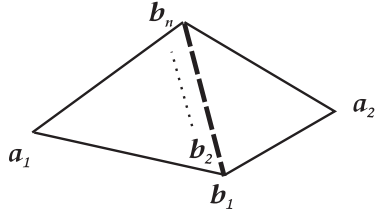}
\caption{Higher order $\mu$-invariants $\mu_j(a_1, \mathfrak{b},
a_2)$.} \label{fig:higher}
\end{figure}
\begin{defin}\label{def:gen-mu} The numbers $\mu_i = \mu_i(a_1, \mathfrak{b},
a_2)$, that appear in the equation {\rm (\ref{eqn:Oda})}, are
referred to as the $\mu$-invariants associated to the collection
$\mathcal{A} =\{a_1, a_2, b_1,\ldots, b_n\}$.
\end{defin}
The following proposition records some of the first, most basic
properties of the $\mu$-invariants $\mu_i(a_1, \mathfrak{b},
a_2)$.

\begin{prop}\label{prop:mu-general-basic}
Given a sequence $\mathfrak{b} = (b_1,\ldots,b_n)$ and a
permutation $\pi$ of integers $1,2,\ldots,n$ let
$\mathfrak{b}^\pi$ be the sequence $(b_{\pi(1)},\ldots
,b_{\pi(n)})$. Suppose that both $(a_1, \mathfrak{b})$ and $(a_2,
\mathfrak{b})$ are ordered bases of the lattice
$\mathbb{Z}^{n+1}$. Then,
\begin{itemize}
 \item[{\rm (a)}] $\mu_i(a_2, \mathfrak{b}, a_1) = -\mu_i(a_1, \mathfrak{b},
a_2)$;
 \item[{\rm (b)}] $\mu_i(a_1, \mathfrak{b}^\pi, a_2) =
 {\rm sgn}(\pi)\mu_i(a_1, \mathfrak{b}, a_2)$.
\end{itemize}
\end{prop}

\medskip\noindent
{\bf Proof:} Both equalities follow from the uniqueness of the
representation (\ref{eqn:Oda}).

\medskip
The $\mu$-invariants are quite natural. They implicitly appear in
many constructions involving lattice polytopes and fans, see for
example the `Oda's criterion' (Theorem~4.12 in
\cite[Chapter~V]{E}). Together with the $\nu$-invariant $\nu(a,
\mathfrak{b})= [a; \mathfrak{b}]$ they are the obvious choice for
the `unimodular discrete curvature invariants' relevant for the
hypothetical generalizations of Theorem~\ref{thm:glavna_teorema}.

This may serve as a justification for a deeper study of
$\mu$-invariants and as an example we prove a generalization of
Proposition~\ref{prop:Jacobi} from Section~\ref{sec:exchange}.

\begin{prop}\label{prop:gen-Jacobi}{\rm (Jacobi identity)}
Suppose that $a_1,a_2,a_3\in \mathbb{Z}^n$. Let
$\mathfrak{b}=(b_1,b_2,\ldots, b_n)$ be a sequence of integer
vectors such that all three collections $(a_1,\mathfrak{b}),
(a_2,\mathfrak{b}), (a_3,\mathfrak{b})$ are bases of
$\mathbb{Z}^n$. Then for each $j\in \{1,\ldots, n\}$
\begin{equation}\label{eqn:gen-Jacobi}
\mu_j(a_1, \mathfrak{b}, a_2) + \mu_j(a_2, \mathfrak{b}, a_3) +
\mu_j(a_3, \mathfrak{b}, a_1) = 0.
\end{equation}
\end{prop}

\medskip\noindent
{\bf Proof:} By adding up the three equations (\ref{eqn:Oda}),
associated to the pairs $(a_1,a_2), (a_2, a_3), (a_3, a_1)$, one
deduces the equality
\begin{equation}
\sum_{j=1}^n(\mu_j(a_1, \mathfrak{b}, a_2) + \mu_j(a_2,
\mathfrak{b}, a_3) + \mu_j(a_3, \mathfrak{b}, a_1))b_j  = 0
\end{equation}
and the result follows from the linear independence of vectors
$b_1, \ldots, b_n$. \hfill $\square$

\medskip\noindent
{\bf Problem~2:} Explore the higher dimensional analogues
(generalizations) of results from Sections~\ref{sec:mu} and
\ref{sec:additivity}.

\medskip
It is not difficult to prove that a unimodular sequence $L =
u_1u_2\ldots u_n$ is completely determined by the initial elements
$u_1$ and $u_2$ and the `discrete curvature invariants'
$\nu(u_i,u_{i+1}), \mu(u_i, u_{i+1}, u_{i+2})$. Moreover, one can
more or less freely prescribe in advance these numbers, as made
precise by the following proposition.

\begin{prop}\label{prop:exists-unique}
 Let $(u_1,u_2)$ be an ordered basis of the lattice
$\mathbb{Z}^2$. Suppose that
\[
\nu_{12}, \nu_{23}\ldots \nu_{n-1,n} \quad \mbox{\rm and} \quad
\mu_2, \mu_3,\ldots , \mu_{n-1}
\]
are sequences of integers such that $\nu_{12}={\rm det}(u_1,u_2)$,
$\nu_{i,i+1}\in \{-1,+1\}$ for each $i=1,\ldots, n-1$ and
$\mu_i\in \mathbb{Z}$ for each $i=2,\ldots,n-1$. Then there exists
a unique unimodular sequence $u_1u_2\ldots u_n$, which extends the
basis $(u_1,u_2)$, such that $\nu_{i,i+1}=\nu(u_i,u_{i+1})={\rm
det}(u_i,u_{i+1})$ for each $i=1,\ldots, n-1$ and $\mu_j =
\mu(u_{j-1},u_j, u_{j+1})$ for each $j=2,\ldots,n-1$.
\end{prop}

\medskip\noindent
{\bf Proof:} Suppose that $u_1,\ldots, u_{k}$ are already
constructed for some $k$, where $2\leq k\leq n-1$. The equation
(\ref{eqn:mu-invarijanta-definicija}), that defines the
$\mu$-invariant, rewritten as follows
\begin{equation}\label{eqn:uni-konstrukcija}
\nu_{k-1,k}u_{k-1} + \nu_{k,k+1}u_{k+1} + \mu_ku_k = 0,
\end{equation}
uniquely determines the vector $u_{k+1}$.

This observation shows how the sequence $(u_i)_1^n$ can be
recursively constructed and it remains to check that
$\nu(u_k,u_{k+1}) = \nu_{k,k+1}$ for each $k=1,\ldots, n-1$ and
$\mu(u_{k-1},u_k, u_{k+1})=\mu_k$ for each $k=2,\ldots,n-1$.

It follows from the equality (\ref{eqn:uni-konstrukcija}) that
\[
\nu_{k-1,k}\nu(u_{k-1},u_k) = \nu_{k,k+1}\nu(u_k,u_{k+1})
\]
and if we already know that $\nu_{k-1,k}=\nu(u_{k-1},u_k)$ we
obtain $\nu_{k,k+1}=\nu(u_{k},u_{k+1})$ as an immediate
consequence. From here and the defining equation
(\ref{eqn:uni-konstrukcija}) we also deduce the equality $\mu_k =
\mu(u_{k-1},u_k,u_{k+1})$. \hfill $\square$

\medskip\noindent
{\bf Problem~3:} Explore the higher dimensional generalizations of
Proposition~\ref{prop:exists-unique}.

\medskip\noindent
{\bf Acknowledgements:} It is a pleasure to acknowledge the
comments and suggestions of the referees which considerably
improved the presentation and organization of the paper.

\end{document}